\documentclass[11pt]{article}
\usepackage{graphicx}
\usepackage{amsfonts}
\usepackage[usenames]{color}

\begin{document}

\title{On the performance of algorithms for the minimization of $\ell_1$-penalized functionals}
\author{Ignace Loris\\Mathematics Department, Vrije Universiteit Brussel\\
Pleinlaan 2, 1050 Brussel, Belgium}
\date{12/12/2008}

\maketitle

\begin{abstract}
The problem of assessing the performance of algorithms used for
the minimization of an $\ell_1$-penalized least-squares
functional, for a range of penalty parameters, is investigated.
A criterion that uses the idea of `approximation isochrones' is
introduced. Five different iterative minimization algorithms
are tested and compared, as well as two warm-start strategies.
Both well-conditioned and ill-conditioned problems are used in
the comparison, and the contrast between these two categories
is highlighted.
\end{abstract}

\section{Introduction}

In recent years, applications of sparse methods in signal
analysis and inverse problems have received a great deal of
attention. The term \emph{compressed sensing} is used to
describe the ability to reconstruct a \emph{sparse} signal or
object from far fewer linear measurements than would be needed
traditionally \cite{Donoho2006}.

A promising approach, applicable to the regularization of
linear inverse problems, consists of using a sparsity-promoting
penalization. A particularly popular penalization of this type
is the $\ell_1$ norm of the object in the basis or frame in
which the object is assumed to be sparse. In
\cite{Daubechies.Defrise.ea2004} it was shown that adding an
$\ell_1$ norm penalization to a least squares functional (see
expression (\ref{l1functional}) below) regularizes ill-posed
linear inverse problems . The minimizer of this functional has
many components exactly equal to zero. Furthermore, the
iterative soft-thresholding algorithm (IST) was shown to
converge in norm to the minimizer of this functional (earlier
work on $\ell_1$ penalties is in
\cite{Santosa.Symes1986,Tibshirani1996}).

It has also been noted that the convergence of the IST
algorithm can be rather slow in cases of practical importance,
and research interest in speeding up the convergence or
developing alternative algorithms is growing
\cite{DaFoL2008,Kim.Koh.ea2007,Figueiredo.Nowak.ea2008,Hale.Yin.ea2007,Beck.Teboulle2008}.
There are already several different algorithms for the
minimization of an $\ell_1$ penalized functional. Therefore, it
is necessary to discuss robust ways of evaluating and comparing
the performance of these competing methods. The aim of this
manuscript is to propose a procedure that assesses the
strengths and weaknesses of these minimization algorithms for a
range of penalty parameters.

Often, authors compare algorithms only for a single value of
the penalty parameter and may thereby fail to deliver a
complete picture of the convergence speed of the algorithms.
For the reader, it is difficult to know if the parameter has
been tuned to favor one or the other method.  Another issue
plagueing the comparison of different minimization algorithms
for problem (\ref{minimizer}) below is the confusion that is
made with sparse recovery. Finding a sparse solution of a
linear equation and finding the minimizer of (\ref{minimizer})
are closely related but (importantly) \emph{different}
problems: The latter will be sparse (the higher the penalty
parameter, the sparser), but the former does not necessarily
minimize the $\ell_1$ penalized functional for any value of the
penalty parameter. Contrary to many discussions in the
literature, we will look at the minimization problem
(\ref{minimizer}) independently of its connection to sparse
recovery and compressed sensing.

The central theme of this note is the introduction of the
concept of \emph{approximation isochrone}, and the illustration
of its use in the comparison of different minimization
algorithms. It proves to be an effective tool in revealing when
algorithms do well and under which circumstances they fail. As
an illustration, we compare five different iterative
algorithms. For this we use a strongly ill-conditioned linear
inverse problem that finds its origin in a problem of seismic
tomography \cite{Loris.Nolet.ea2007}, a Gaussian random matrix
(i.e. the matrix elements are taken from a normal distribution)
as well as two additional synthetic matrices. In the existing
literature, most tests are done using only a matrix of random
numbers, but we believe that it is very important to also
consider other matrices. Actual inverse problems may depend on
an operator with a less well-behaved spectrum or with other
properties that could make the minimization more difficult.
Such tests are usually not available. Here we compare four
operators. Among other things, we find that the strongly
singular matrices are more demanding on the algorithms.

We limit ourselves to the case of real signals, and do not
consider complex variables. In this manuscript, the usual
$2$-norm of a vector $x$ is denoted by $\|x\|$ and the $1$-norm
is denoted by $\|x\|_1$. The largest singular value of a matrix
$K$ is denoted by $\|K\|$.

\section{Problem statement}

\label{problemsection}

After the introduction of a suitable basis or frame for the
object and the image space, the minimization problem under
study can be stated in its most basic form, without referring
to any specific physical origin, as the minimization of the
convex functional
\begin{equation}
F_\lambda(x)=\|Kx-y\|^2+2\lambda \|x\|_1\label{l1functional}
\end{equation}
in a real vector space ($x\in \mathbb{R}^p$, $\lambda\geq 0$),
for a fixed linear operator $K\in \mathbb{R}^{m\times p}$ and
data $y\in \mathbb{R}^m$.   In the present analysis we will
assume the linear operator $K$ and the data $y$ are such that
the minimizer of (\ref{l1functional}) is unique. This is a
reasonable assumption as one imposes penalty terms, typically,
to make the solution to an inverse problem unique.

We set
\begin{equation}
\bar x(\lambda)=\arg\min_x\|Kx-y\|^2+2\lambda
\|x\|_1.\label{minimizer}
\end{equation}
The penalty parameter $\lambda$ is positive; in applications,
it has to be chosen depending on the context. Problem
(\ref{minimizer}) is equivalent to the constrained minimization
problem:
\begin{equation}
\tilde x(\rho)=\arg\min_{\|x\|_1\leq \rho}\|Kx-y\|^2 \qquad
\label{constrmin}
\end{equation}
with an implicit relationship between $\rho$ and $\lambda$:
$\rho=\|\bar x(\lambda)\|_1$. It follows from the equations
(\ref{kkt}) below that the inverse relationship is:
$\lambda=\max_i |(K^T(y-K \tilde x(\rho)))_i|$. Under these
conditions one has that: $\bar x(\lambda)=\tilde x(\rho)$. One
also has that $\bar x(\lambda)=0$ for all $\lambda\geq
\lambda_\mathrm{max}\equiv\max_i|(K^Ty)_i|$.

\subsection{Direct method}

An important thing to note is that the minimizer $\bar x$ (and
thus also $\tilde x$) can in principle be found in a finite
number of steps using the homotopy/LARS method
\cite{Osborne.Presnell.ea2000,Efron.Hastie.ea2004}. This direct
method starts from the variational equations which describe the
minimizer $\bar x$:
\begin{equation}
\begin{array}{lclcl}
(K^T(y-K\bar x))_i&=&\lambda\; \mathrm{sgn}(\bar x_i) & \mathrm{if} & \bar x_i\neq 0 \\
|(K^T(y-K\bar x))_i|&\leq&\lambda & \mathrm{if} & \bar x_i= 0
\end{array}\label{kkt}
\end{equation}
Because of the piece-wise linearity of the equations
(\ref{kkt}), it is possible to construct $\bar
x(\lambda_\mathrm{stop})$ by letting $\lambda$ in (\ref{kkt})
descend from $\lambda_\mathrm{max}$ to $\lambda_\mathrm{stop}$,
and by solving a linear system at every value of $\lambda$
where a component in $\bar x(\lambda)$ goes from zero to
nonzero or, exceptionally, from nonzero to zero. The first such
point occurs at $\lambda=\lambda_\mathrm{max}$. The linear
systems that have to be solved at each of these breakpoints are
`small', starting from $1\times 1$ and ending with $s\times s$,
where $s$ is the number of nonzero components in $\bar
x(\lambda_\mathrm{stop})$. Such a method thus constructs $\bar
x(\lambda)$ exactly for all
$\lambda_\mathrm{max}\geq\lambda\geq \lambda_\mathrm{stop}$, or
equivalently all $\bar x$ with $0\leq\|\bar x(\lambda)\|_1\leq
\|\bar x(\lambda_\mathrm{stop})\|_1$. It also follows that
$\bar x(\lambda)$ is a piecewise linear function of $\lambda$.

\begin{figure}
\centering\resizebox{\textwidth}{!}{\includegraphics{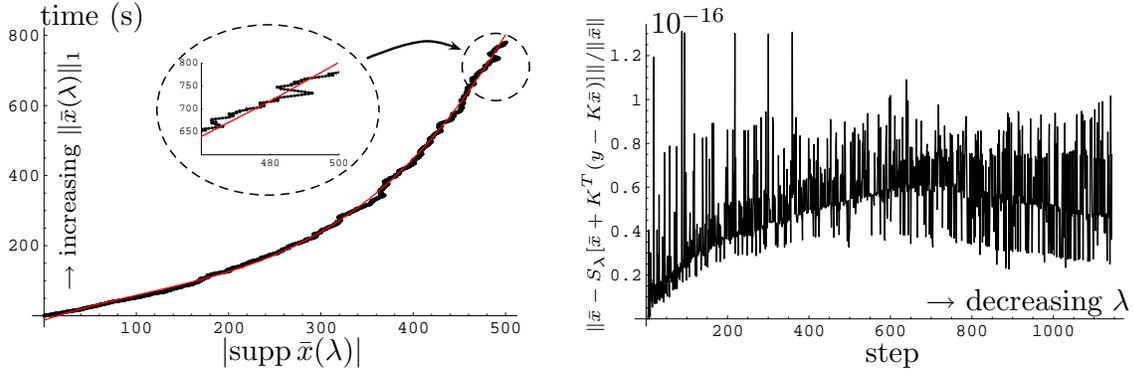}}
\caption{Left: Time complexity of the direct algorithm
mentioned in section \ref{problemsection}. Horizontal axis:
number of nonzero components in the minimizer $\bar x$.
Vertical axis: time needed by the direct algorithm to calculate
this minimizer (seconds). The continuous line represents a
cubic fit. The support of $\bar x(\lambda)$ may sometimes
decrease for increasing values of $\|\bar x(\lambda)\|_1$. This
phenomenon is clearly visible in the zoomed area. Right: The
relative error, $\|\bar x-S_\lambda(\bar x+K^T(y-K\bar
x))\|/\|\bar x\|$, of the resulting minimizer at each step
(there are more steps than nonzero components
because of the phenomenon in the zoomed area in the left hand
side plot). Due to floating point arithmetic, this is not
exactly zero. These two pictures were made by letting $\lambda$
decrease from $\lambda_\mathrm{max}$ to
$\lambda_\mathrm{stop}=\lambda_\mathrm{max}/2^{15.5654}$. The operator used is $K^{(1)}$.}\label{complexityfig}
\end{figure}

Implementations of this direct algorithm exist
\cite{IMM2005-03897,Donoho.Stodden.ea2007,Loris2007}, and
exhibit a time complexity that is approximately cubic in $s$
(the number of nonzero components of $\bar x(\lambda)$). If one
is interested in sparse recovery, this is not necessarily a
problem as time complexity is linear for small $s$. In fact,
the algorithm can be quite fast, certainly if one weighs in the
fact that the exact minimizer (up to computer round-off) is
obtained. A plot of the time complexity as a function of the
number of nonzero components in $\bar x(\lambda)$ is given in
figure \ref{complexityfig} (left hand side) for an example
operator $K^{(1)}$ and data $y$ (see start of section
\ref{comparisonsection} for a description of $K^{(1)}$). The
plot shows that the direct algorithm is useable in practice for
about $s\leq 10^3$. This graph also illustrates the fact that
the size of the support of $\bar x(\lambda)$ does not
necessarily grow monotonically with decreasing penalty
parameter (this depends on the operator and data).

It follows immediately from equation (\ref{kkt}) that the minimizer
$\bar x(\lambda)$ satisfies the fixed point equation
\begin{equation}
\bar x=S_\lambda[\bar x+K^T(y-K\bar x)]\label{fixedpoint}
\end{equation}
where $S_\lambda$ is the well-known soft-thresholding operator
applied component-wise:
\begin{equation}
S_\lambda(u)=\left\{\begin{array}{llccr}
u-\lambda & \qquad\qquad & u&\geq&\lambda\\
0 & & |u|&\leq&\lambda\\
u+\lambda & & u&\leq& -\lambda
\end{array}\right.
\end{equation}
In real life we have to take into account that computers work
with floating point, inexact, arithmetic. The definition of an
exact solution, in this case, is that $\bar x(\lambda)$
satisfies the fixed-point equation (\ref{fixedpoint}) up to
computer precision:
\begin{equation}
\|\bar x-S_\lambda[\bar x+K^T(y-K\bar x)]\|/\|\bar x\| \approx10^{-16}.\label{inexactfixedpoint}
\end{equation}
The direct algorithm mentioned before can be implemented using
floating point arithmetic
(\cite{IMM2005-03897,Donoho.Stodden.ea2007} do this). The
implementation \cite{Loris2007} can handle both exact
arithmetic (with integer, rational numbers) and floating point
arithmetic. An example of the errors made by the direct method
is illustrated figure \ref{complexityfig} (right). We will
still use the term `exact solution' as long as condition
(\ref{inexactfixedpoint}) is satisfied.

\subsection{Iterative algorithms}

There exist several iterative methods that can be used for the
minimization problems (\ref{minimizer}) or (\ref{constrmin}):
\begin{enumerate}
\item The iterative soft-thresholding (IST) algorithm
    already mentioned in the introduction can be written
    as:
\begin{equation}
x^{(n+1)}=S_\lambda[x^{(n)}+K^T(y-Kx^{(n)})],\qquad
x^{(0)}=0.\label{tlw}
\end{equation}
Under the condition $\|K\|<1$ the limit of this sequence coincides
with the minimizer $\bar x(\lambda)$ and $F_\lambda(x^{(n)})$
decreases monotonically as a function of $n$
\cite{Daubechies.Defrise.ea2004}. For $\|K\|<\sqrt{2}$, there is
still (weak) convergence \cite[Corollary 5.9]{Combettes.Wajs2005},
but the functional (\ref{l1functional}) is no longer guaranteed to
decrease at every step. This algorithm is probably the easiest to
implement.\label{tlwalg}

\item a projected steepest descent method \cite{DaFoL2008}
    (and related \cite[expression (59)]{Combettes1997}):
\begin{equation}
x^{(n+1)}=P_\rho[x^{(n)}+\beta^{(n)}\,r^{(n)}],\qquad
x^{(0)}=0,
\end{equation}
with $r^{(n)}=K^T(y-Kx^{(n)})$ and
$\beta^{(n)}=\|r^{(n)}\|^2/\|Kr^{(n)}\|^2$. $P_\rho(\cdot)$
denotes the orthogonal projection onto an $\ell_1$ ball of
radius $\rho$, and can be implemented efficiently by
soft-thresholding with an appropriate variable
threshold.\label{psdalg}

\item the `GPSR-algorithm' (gradient projection for sparse
    reconstruction), another iterative projection method,
    in the auxiliary variables $u,v\geq 0$ with $x=u-v$
    \cite{Figueiredo.Nowak.ea2008}.\label{gpsralg}

\item the `$\ell_1$-ls-algorithm', an interior point method
    using preconditioned conjugate gradient substeps (this
    method solves a linear system in each outer iteration
    step) \cite{Kim.Koh.ea2007}.\label{l1lsalg}

\item `FISTA' (fast iterative soft-thresholding algorithm)
    is a variation of the iterative soft-thresholding
    algorithm. Define the (non-linear) operator $T$ by
    $T(x)=S_\lambda[x+K^T(y-Kx)]$. Then the FISTA algorithm
    is:
\begin{equation}
x^{(n+1)}=T\left(x^{(n)}+\frac{t^{(n)}-1}{t^{(n+1)}} \left(x^{(n)}-x^{(n-1)}\right)\right)\qquad
x^{(1)}=0,
\label{fista}
\end{equation}
where $t^{(n+1)}=\frac{1+\sqrt{1+4(t^{(n)})^2}}{2}$ and
$t^{(1)}=1$. It has virtually the same complexity as
algorithm \ref{tlwalg}, but can be shown to have better
convergence properties
\cite{Beck.Teboulle2008}.\label{fistaalg}
\end{enumerate}

\subsection{Warm-start strategies}

There also exist so-called \emph{warm-start} strategies. These
methods start from $\bar x(\lambda_0=\lambda_\mathrm{max})=0$
and try to approximate $\bar x(\lambda_k)$  for $k:0,\ldots,N$
by starting from an approximation of $\bar x(\lambda_{k-1})$
already obtained in the previous step instead of always
restarting from $0$. They can be used for finding an
approximation of a whole range of minimizers $\bar
x(\lambda_k)$ for a set of penalty parameters
$\lambda_\mathrm{max}=\lambda_0>\lambda_1>\lambda_2>\ldots>\lambda_N=\lambda_\mathrm{stop}$
or, equivalently, for a set of $\ell_1$-radii
$0=\rho_0<\rho_1<\rho_2<\ldots<\rho_N=\rho_\mathrm{stop}$. Two
examples of such methods are:
\begin{enumerate}
\item[(A)] `fixed-point continuation' method
    \cite{Hale.Yin.ea2007}:
\begin{equation}
x^{(n+1)}=S_{\lambda_{n+1}}[x^{(n)}+K^T(y-Kx^{(n)})]
\end{equation}
with $\lambda_0=\lambda_\mathrm{max}$ and
$\lambda_{n+1}=\alpha\lambda_{n}$ and $\alpha<1$ such that
$\lambda_N=\lambda_\mathrm{stop}$ (after a pre-determined
number $N$ of steps). In other words, the threshold is
decreased geometrically instead of being fixed as in the
IST method and $x^{(n)}$ is interpreted as an approximation
of $\bar x(\lambda_n)$.

\item[(B)] adaptive steepest descent \cite{DaFoL2008}:
\begin{equation}
x^{(n+1)}=P_{\rho_{n+1}}[x^{(n)}+\beta^{(n)}\,K^T(y-Kx^{(n)})],\qquad
x^{(0)}=0,
\end{equation}
with $r^{(n)}=K^T(y-Kx^{(n)})$,
$\beta^{(n)}=\|r^{(n)}\|^2/\|Kr^{(n)}\|^2$, and
$\rho_{n+1}=(n+1) \rho_\mathrm{stop}/N$. Here the radius
$\rho_n$ increases arithmetically instead of being fixed as
in algorithm \ref{psdalg} and $x^{(n)}$ is interpreted as
an approximation of $\tilde x(\rho_n)$.
\end{enumerate}
Such algorithms have the advantage of providing an
approximation of the Pareto curve (a plot of $\|K\bar
x(\lambda)-y\|^2$ vs. $\|\bar x(\lambda)\|_1$, also known as
trade-off curve) as they go, instead of just calculating the
minimizer corresponding to a fixed penalty parameter. It is
useful for determining a suitable value of the penalty
parameter $\lambda$ in applications.

\section{Approximation isochrones}

\label{isosection}

In this section, we discuss the problem of assessing the speed
of convergence of a given minimization algorithm.

The minimization problem (\ref{l1functional}) is often used in
compressed sensing. In this context, an iterative minimization
algorithm may be tested as follows: one chooses a (random)
sparse input vector $x^\mathrm{input}$, calculates the image
under the linear operator and adds noise
$y=Kx^\mathrm{input}+n$. One then uses the algorithm in
question to try to reconstruct the input vector
$x^\mathrm{input}$ as a minimizer of (\ref{minimizer}),
choosing $\lambda$ in such a way that the resulting $x^{(N)}$
best corresponds to the input vector $x^\mathrm{input}$. This
procedure in not useful in our case because it does not compare
the iterates $x^{(n)}$ ($n:0\ldots N$) with the actual
minimizer $\bar x$ of functional (\ref{minimizer}). Even though
it is sparse, the input vector $x^\mathrm{input}$ most likely
does not satisfy equations of type (\ref{kkt}), and hence does
not constitute an actual minimizer of (\ref{l1functional}).
Such a type of evaluation is e.g. done
in \cite[section IV.A]{Figueiredo.Nowak.ea2008}.\\
In this note we are interested in describing how well an
algorithm does in finding the true minimizer of
(\ref{l1functional}), not in how suitable an algorithm may be
for compressed sensing applications. We consider sparse
recovery and $\ell_1$-penalized functional minimization as two
separate issues. Here, we want to focus on the latter.

Another unsatisfactory method of evaluating the convergence of
a minimization algorithm is to look at the behavior of the
functional $F_\lambda(x^{(n)})$ as a function of $n$. For small
penalties, it is quite possible that the minimum is almost
reached by a vector $x^{(N)}$ that is still quite far from the
real minimizer $\bar x$.

Suppose one has developed an iterative algorithm for the
minimization of the $\ell_1$-penalized functional
(\ref{l1functional}), i.e. a method for the computation of
$\bar x(\lambda)$ in expression (\ref{minimizer}). As we are
interested in evaluating an algorithm's capabilities of
minimizing the functional (\ref{l1functional}), it is
reasonable that one would compare the iterates $x^{(n)}$ with
the exact minimizer $\bar x(\lambda)$. I.e. evaluation of the
convergence speed should look at the quantity $\|x^{(n)}-\bar
x\|$ as a function of time. A direct procedure exists for
calculating $\bar x(\lambda)$ and thus it is quite
straightforward to make such an analysis for a whole range of
values of the penalty parameter $\lambda$ (as long as the
support of $\bar x(\lambda)$ is not excessively large). We will
see that the weaknesses of the iterative algorithms are already
observable for penalty parameters corresponding to quite sparse
$\bar x(\lambda)$, i.e. for $\bar x$ that are still relatively
easy to compute with the direct method.

In doing so, one has three parameters that need to be included
in a graphical representation: the relative reconstruction
error $e=\|x^{(n)}-\bar x\|/\|\bar x\|$, the penalty parameter
$\lambda$ and also the time $t$ needed to obtain the
approximation $x^{(n)}$. Making a 3D plot is not a good option
because, when printed on paper or displayed on screen, it is
difficult to accurately interpret the shape of the surface. It
is therefore advantageous to try to make a more condensed
representation of the algorithm's outcome. One particularly
revealing possibility we suggest, is to use the
$\lambda$-$e$-plane to plot the \emph{isochrones} of the
algorithm.

For a fixed amount of computer time $t$, these isochrones trace
out the approximation accuracy $e$ that can be obtained for
varying values of the penalty parameter $\lambda$. There are
two distinct advantages in doing so. Firstly, it becomes
immediately clear for which range of $\lambda$ the algorithm in
question converges quickly and, by labeling the isochrones, in
which time frame. Secondly, it is clear where the algorithm has
trouble approaching the real minimizer: this is characterized
by isochrones that are very close to each other and away from
the $e=0$ line. Hence a qualitative evaluation of the
convergence properties can be made by noticing where the
isochrones are broadly and uniformly spaced (good convergence
properties), and where the isochrones seem to stick together
(slow convergence). Quantitatively, one immediately sees what
the remaining relative error is.

Another advantage of this representation is that (in some weak
sense) the isochrones do not depend on the computer used: if
one uses a faster/slower computer to make the calculations, the
isochrones `do not change place' in the sense that only their
labels change. In other words, this way of representing the
convergence speed of the algorithm accurately depicts the
region where convergence is fast or slow, independently of the
computer used.

\begin{figure}
\centering\resizebox{\textwidth}{!}{\includegraphics{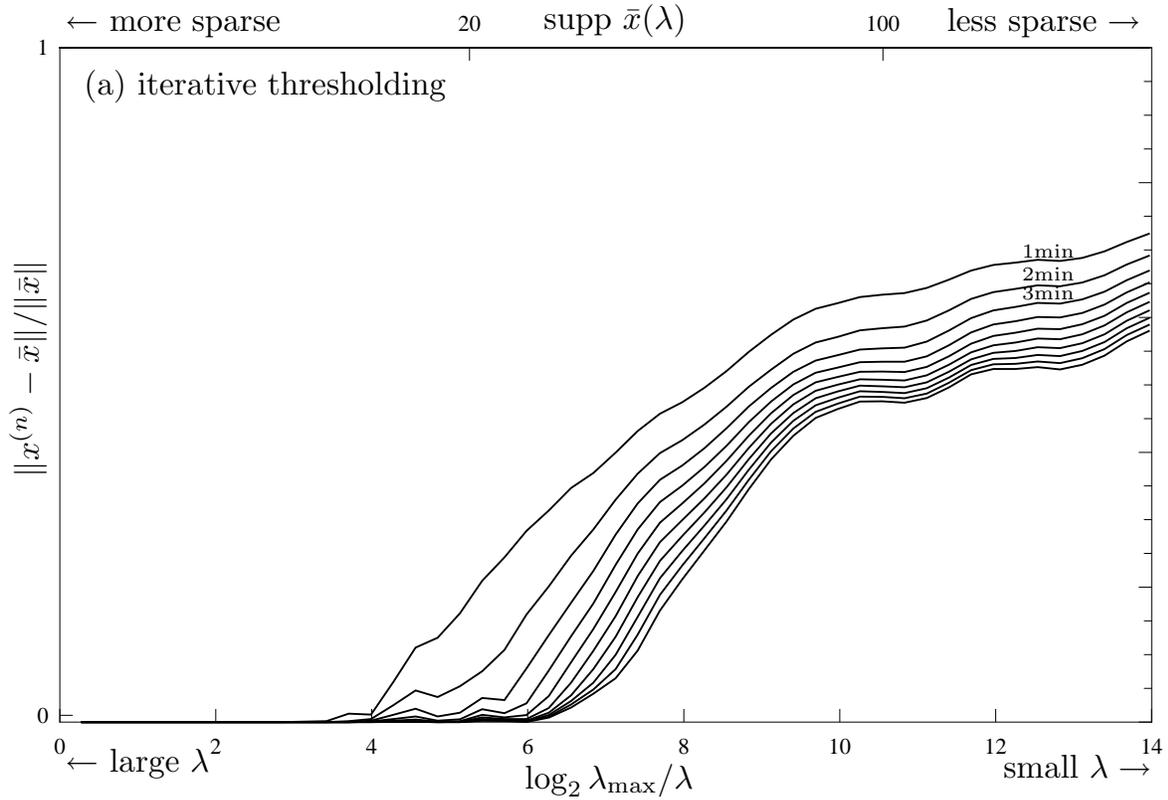}}
\caption{This figure displays the approximation isochrones
($t=1,2,\ldots,10$ minutes) for the iterative thresholding
algorithm (\ref{tlw}) with operator $K^{(1)}$. The vertical axis is $e=\|x^{(n)}-\bar
x\|/\|\bar x\|$. The bottom horizontal axis represents
$\log_2\lambda_\mathrm{max}/\lambda$ (i.e. small $\lambda$ to
the right). The top horizontal axis is used to indicate the number
of nonzero components in the corresponding minimizer (nonlinear
scale). Convergence is satisfactory for
$\lambda>\lambda_\mathrm{max}/2^6$ but very slow for
$\lambda<\lambda_\mathrm{max}/2^8$}\label{singleisochronefig}
\end{figure}

In figure \ref{singleisochronefig}, an example of such a plot
is given. The operator is again $K^{(1)}\in
\mathbb{R}^{1848\times 8192}$ (more details are given at the
start of section \ref{comparisonsection}). The algorithm
assessed in this plot is the iterative thresholding algorithm
(a). Clearly one sees that the iterative thresholding does well
for $\lambda\geq \lambda_\mathrm{max}/2^{6}$, but comes into
trouble for $\lambda\leq \lambda_\mathrm{max}/2^{8}$. This
means that the iterative thresholding algorithm (\ref{tlw}) has
trouble finding the minimizer with more than about 50 nonzero
components (out of a possible 8192 degrees of freedom). The
direct method is still practical in this regime as proven by
figure \ref{complexityfig} (left), and it was used to find the
`exact' $\bar x(\lambda)$'s. In fact, here the direct method is
faster than the iterative methods.

\section{Comparison of minimization algorithms}

\label{comparisonsection}

In this section we will compare the six iterative algorithms
(a)--(e), and the two warm-start algorithms (A)--(B) mentioned
in section \ref{problemsection}. We will use four qualitatively
different operators $K$ for making this comparison.

\begin{figure}
\centering\resizebox{\textwidth}{!}{\includegraphics{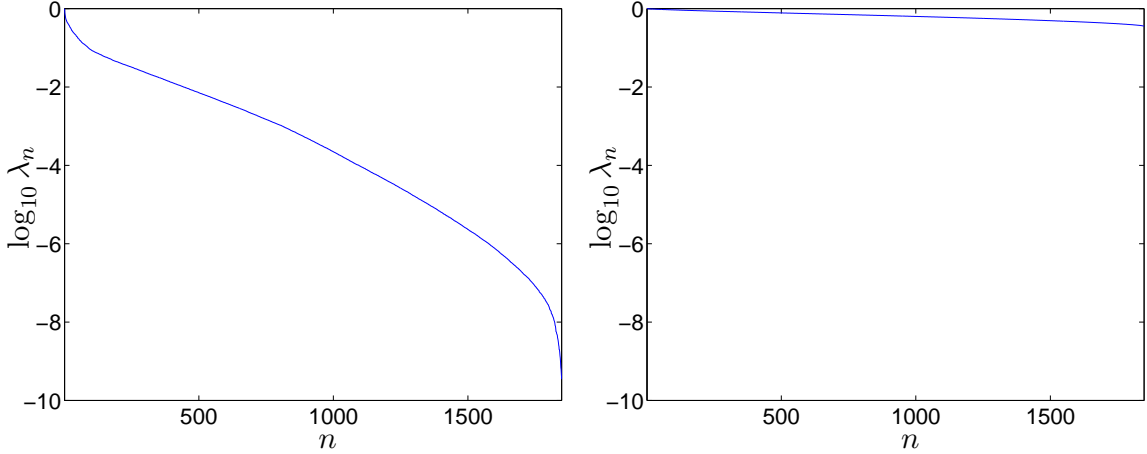}}
\caption{Left: The singular values $\lambda_n$ ($n:1,\ldots, 1848$) of the operators $K^{(1)}$
and $K^{(4)}$. Right: The singular values $\lambda_n$ ($n:1,\ldots, 1848$) of the operators
$K^{(2)}$ and $K^{(3)}$. The Gaussian random matrix $K^{(2)}$ is much better conditioned than operator $K^{(1)}$.}\label{spectrumpic}
\end{figure}

Firstly, we will use an ill-conditioned matrix $K^{(1)}$
stemming from a geo-science inverse problem. It was already
used in figures \ref{complexityfig} and
\ref{singleisochronefig}. It contains $1848$ 2-D integration
kernels discretized on a $64\times 64$ grid, and expanded in a
($2\times$ redundant) wavelet frame. Hence this matrix has
$1848$ rows and $8192$ ($=2\times 64^2$) columns. The spectrum
is pictured in figure~\ref{spectrumpic}. Clearly it is severely
ill-conditioned.

The matrix $K^{(2)}$ is of the same size, but contains random
numbers taken from a Gaussian distribution. Its spectrum is
also in figure~\ref{spectrumpic} and has a much better
condition number (ratio of largest singular value to smallest
nonzero singular value). This type of operator is often used in
the evaluation of algorithms for minimization of an
$\ell_1$-penalized functional.

As we shall see, the different minimization algorithms will not
perform equally well for these two operators. In order to
further discuss the influence of both the spectrum and
orientation of the null space of the operator on the
algorithms' behavior, we will also use two other operators.
$K^{(3)}$ will have the same well behaved spectrum as
$K^{(2)}$, but an unfavorably oriented null space: the null
space will contain many directions that are almost parallel to
an edge or a face of the $\ell_1$ ball. $K^{(4)}$ will have the
same ill-behaved spectrum as $K^{(1)}$, but it will have the
same singular vectors as the Gaussian matrix $K^{(2)}$.

More precisely, $K^{(3)}$ is constructed artificially from
$K^{(2)}$ by setting columns $4000$ through $8192$ equal to
column $4000$. This creates an operator with a null space that
contains many vectors parallel to a side or edge of the
$\ell_1$ ball. A small perturbation is added in the form of
another random Gaussian matrix to yield the intermediate matrix
$A$. The singular value decomposition is calculated $A=USV^T$
($U^{-1}=U^T$ and $V^{-1}=V^T$). In this decomposition, the
spectrum is replaced by the spectrum of $K^{(2)}$. This then
forms the matrix $K^{(3)}$.

$K^{(4)}$ is constructed by calculating the singular value
decomposition of $K^{(2)}=USV^T$ and replacing the singular
values in $S$ by those of $K^{(1)}$.

In all cases, $K^{(i)}$ is normalized to have its largest
singular value  $\approx 0.999$ (except when studying algorithm
(a') where we use normalization $\approx 0.999\times\sqrt{2}$).

\subsection{A severely ill-conditioned operator}

\label{illcondsection}

\begin{figure}
\centering\resizebox{\textwidth}{!}{\includegraphics{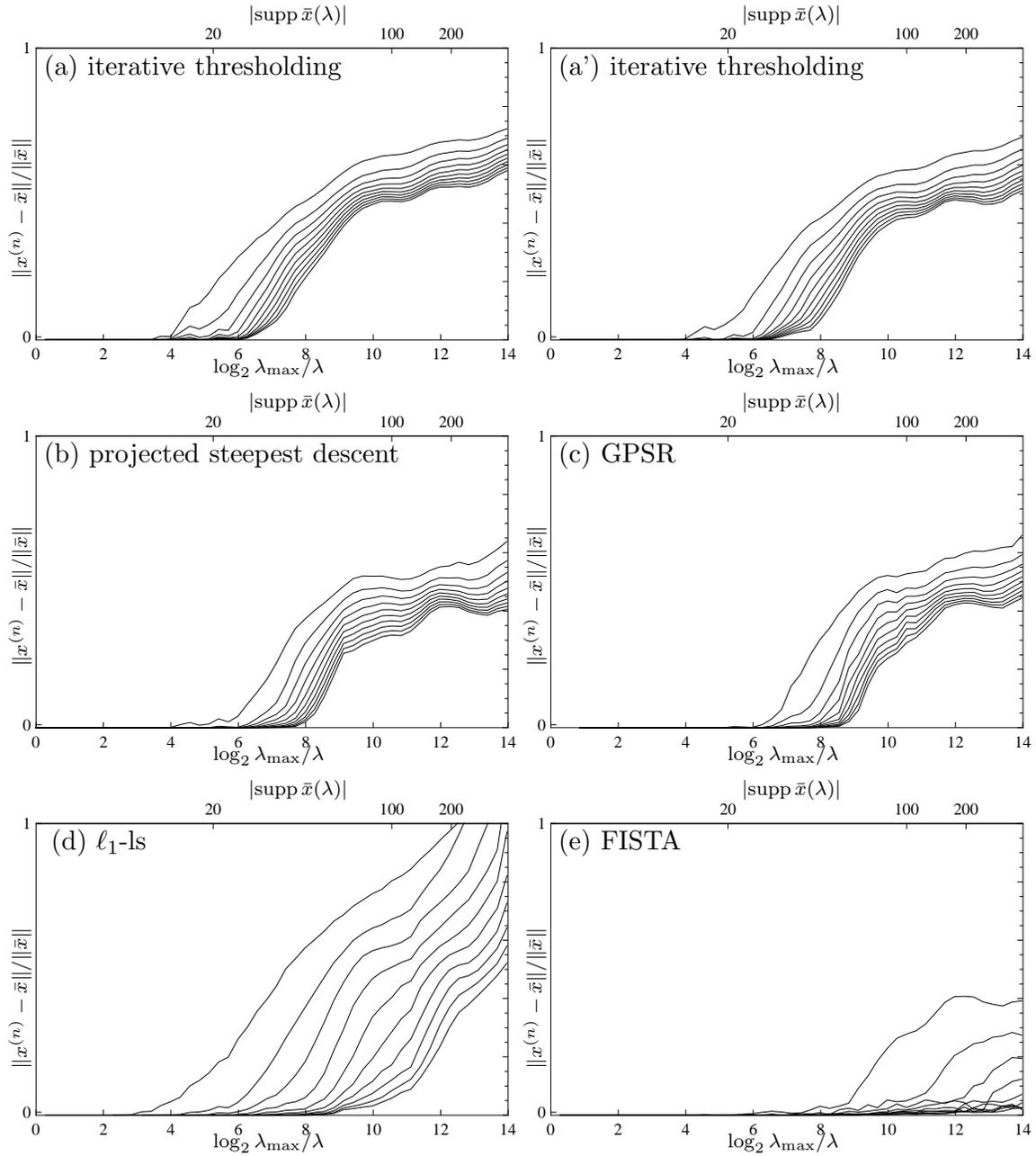}}
\caption{These pictures contain the approximation isochrones
for the algorithms (a)--(e) mentioned in section
\ref{problemsection} for $t=1,2\ldots 10$ minutes. The
horizontal and vertical axis are identical to the ones used in
figure \ref{singleisochronefig}. The operator used is $K^{(1)}$. Clearly, for this example,
methods (a), (a'), (b) and (c) have a lot of difficulty
approaching the minimizer for small values of $\lambda$
(closely spaced isochrones). The $\ell_1$-ls method (d) still
works well, but it is slower for large penalties. The Fista
methods (e) appears to work best for small penalty parameters. The GPSR method works well for relatively large values of $\lambda$.
See text for a discussion.}\label{isochronefig}
\end{figure}

In figure \ref{isochronefig}, we compare the six algorithms
(a)--(e) mentioned in section \ref{problemsection}
 for the same ill-conditioned operator $K^{(1)}\in \mathbb{R}^{1848\times
8192}$ as in figure \ref{singleisochronefig}. We again choose
penalty parameters
$\lambda_\mathrm{max}\geq\lambda\geq\lambda_\mathrm{max}/2^{14}$
and show the isochrones corresponding to $t=1,2,\ldots,10$
minutes. Panel (a) is identical to figure
\ref{singleisochronefig}. All six algorithms do well for large
penalties (left-hand sides of the graphs). For smaller values
of $\lambda$ ($8<\log_2 \lambda_\mathrm{max}/\lambda$) the
isochrones come closer together meaning that convergence
progresses very slowly for algorithms (a), (a'), (b) and (c).
For algorithm (d), the isochrones are still reasonable
uniformly spaced even for smaller values of the penalty
parameter. In this case the projected algorithms do better than
iterative thresholding, but the $\ell_1$-ls algorithm (d) is to
be preferred in case of small penalty parameters. The Fista (e)
algorithm seems to perform best of all for small parameters
$\lambda$.

Apart from the shape of the isochrone curves, it is also
important to appreciate the top horizontal scales of these
plots. The top scale indicate the size of the support of the
corresponding minimizer $\bar x(\lambda)$. We see that all
algorithms have much difficulty in finding minimizers with more
than about $100$ nonzero components. In this range of the
number of nonzero components in $\bar x(\lambda)$, the direct
method is faster for $K^{(1)}$.

\begin{figure}
\centering\resizebox{\textwidth}{!}{\includegraphics{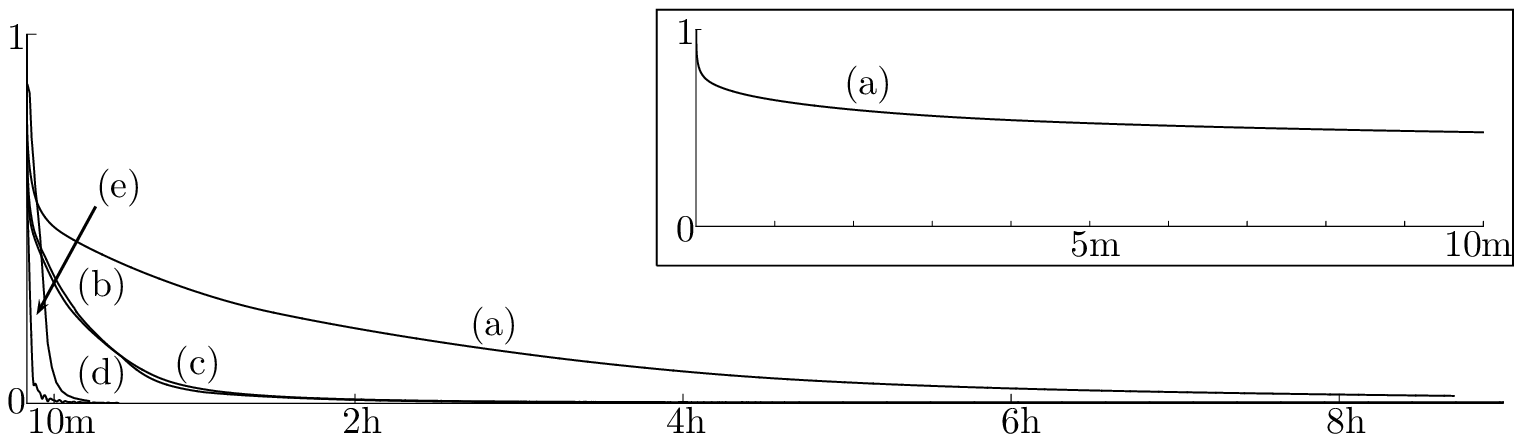}}
\caption{Main graph: Relative error $\|x^{(n)}-\bar x\|/\|\bar
x\|$ as a function of time for one particular value of
$\lambda$. The operator is $K^{(1)}$. Clearly, all algorithms converge to the same limit
$\bar x$, but some algorithms are very slow. Inset: zoom-in on
the first 10 minutes for the iterative thresholding algorithm.
It appears to indicate convergence (to a different minimizer),
but this is deceptive as the main graph shows. Labels (a)--(e)
refer to the algorithms mentioned in section
\ref{problemsection} and figure
\ref{isochronefig}.}\label{singlelambdaplot}
\end{figure}

A skeptic might argue that, in the case of figures
\ref{singleisochronefig} and \ref{isochronefig}, the minimizer
$\bar x(\lambda)$ might not be unique for small values of
$\lambda$, and that this is the reason why the isochrones do
not tend to $e=0$ after about 10 minutes ($\approx 4500$
iterations). This is not the case. If one runs the iterative
methods for a much longer time, one sees that the error $e$
does go to zero. Such a plot is made in figure
\ref{singlelambdaplot} for one choice of $\lambda$
($\lambda\approx\lambda_\mathrm{max}/2^{11.115}$).

What we notice here is that, for iterative soft-thresholding,
$\|x^{(n+1)}-x^{(n)}\|/\|x^{(n)}\|=\|x^{(n)}-S_\lambda[x^{(n)}+K^T(y-K
x^{(n)})\|/\|x^{(n)}\|$ is small ($\approx 10^{-5}$ after
$4500$ iterations in this example). But this does \emph{not}
mean that the algorithm has almost converged (as might be
suggested by figure \ref{singlelambdaplot}-inset); on the
contrary it indicates that the algorithm is progressing only
very slowly for this value of $\lambda$! The difference
$\|x^{(n+1)}-x^{(n)}\|/\|x^{(n)}\|$ should be of the order
$10^{-16}$, for one to be able to conclude convergence, as
already announced in formula (\ref{inexactfixedpoint}).


\subsection{A Gaussian random matrix}

\label{wellcondsection}

In this subsection we choose $K=K^{(2)}$. It is much less
ill-conditioned than the matrix in the previous subsection.

\begin{figure}
\centering\resizebox{\textwidth}{!}{\includegraphics{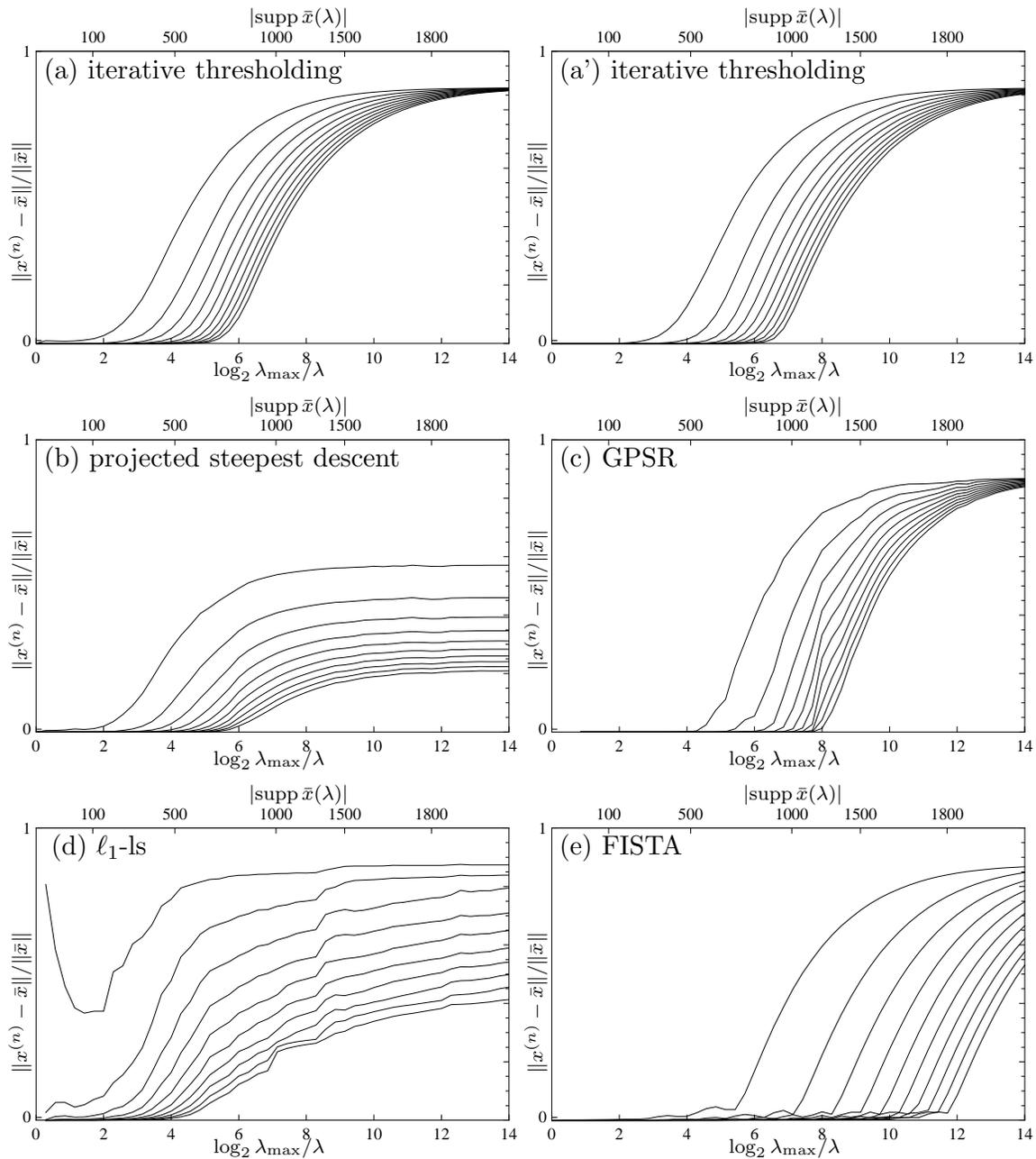}}
\caption{These pictures contain the approximation isochrones
for the algorithms (a)--(e) mentioned in section
\ref{problemsection} for $t=6,12,\ldots, 60$ seconds. The
operator used in this comparison is the Gaussian random matrix $K^{(2)}$.
This matrix is much better conditioned than the matrix $K^{(1)}$ which was used for
figures \ref{isochronefig}. The main differences are faster
convergence and, importantly, more nonzero components are
recoverable (the top scales, indicating the size of the support of the
minimizer, are much larger than in figures \ref{isochronefig}).
}\label{isochronefig3}
\end{figure}

In figure \ref{isochronefig3} we make the same comparison of
the six iterative algorithms (a)--(e) for the operator
$K^{(2)}$. We choose penalty parameters
$\lambda_\mathrm{max}\geq\lambda\geq\lambda_\mathrm{max}/2^{14}$
and show the isochrones corresponding to $t=6,12,\ldots,60$
seconds. I.e. the time scale is $10$ times smaller than for the
ill-conditioned matrix $K^{(1)}$ in the previous section. Again
all the algorithms do reasonably well for large penalty
parameters, but performance diminishes for smaller values of
$\lambda$.

The iterative soft-thresholding method with
$\|K^{(2)}\|\approx\sqrt{2}$ in (a') does slightly better than
iterative soft-thresholding with $\|K^{(2)}\|\approx 1$ in
panel (a). The GPSR method (c) does better than the other
methods for large values of the penalty (up to
$|\mathrm{supp}\, \bar x|\approx 500$), but loses out for
smaller values. The $\ell_1$-ls method (d) does well compared
to the other algorithms as long as the penalty parameter is not
too large. The FISTA method (e) does best except for the very
small value of $\lambda$ (right hand sides of plots) where the
projected steepest descent does better in this time scale.

Apart from the different time scales ($1$ minute vs. $10$
minutes) there is another, probably even more important,
difference between the behavior of the algorithms for $K^{(1)}$
and $K^{(2)}$. In the latter case the size of the support of
the minimizers that are recoverable by the iterative algorithms
range from $0$ to about $1800$ (which is about the maximum for
$1848$ data). This is much more than in the case of the matrix
$K^{(1)}$ in figure \ref{isochronefig} where only minimizers
with about $120$ nonzero coefficients were recoverable.

\subsection{Further examples}

Here we compare the various iterative algorithms for the
operators $K^{(3)}$ and $K^{(4)}$. The former was constructed
such that many elements of its null space are almost parallel
to the edges of the $\ell_1$ ball, in an effort to make the
minimization (\ref{l1functional}) more challenging. The latter
operator has the same singular vectors as the random Gaussian
matrix $K^{(2)}$ but the ill-conditioned spectrum of $K^{(1)}$.
This will then illustrate how the spectrum of $K$ can influence
the algorithms used for solving (\ref{l1functional}).

\begin{figure}
\centering \resizebox{\textwidth}{!}{\includegraphics{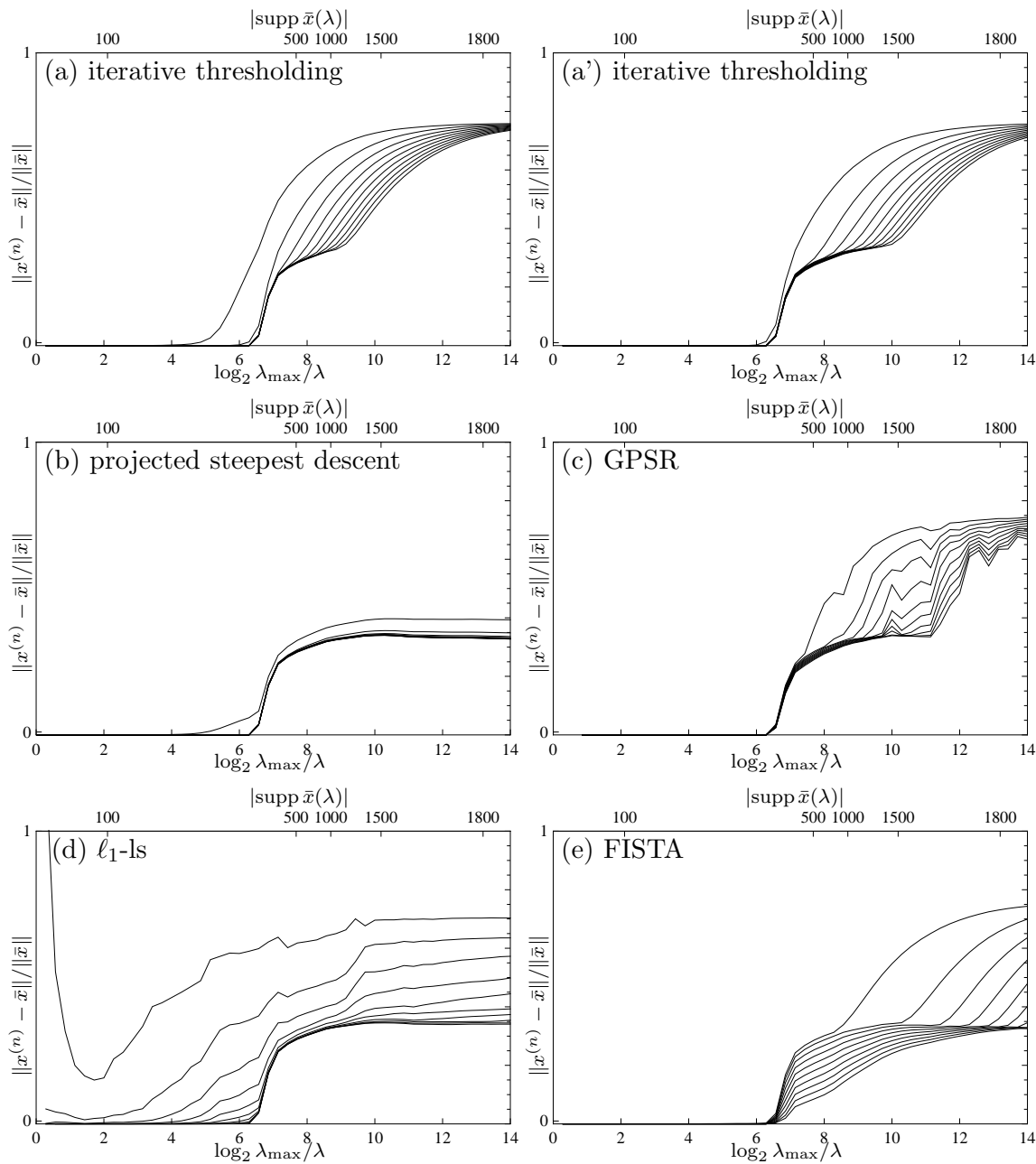}}
\caption{The six algorithms (a)--(e) are compared for the operator $K^{(3)}$.
The time intervals are again $6,12,\ldots,60$s.
For $\lambda<\lambda_\mathrm{max}/2^6$, the algorithms performs
worse than in the case of matrix $K^{(2)}$ in
figure~\ref{isochronefig3}, although $K^{(2)}$ and $K^{(3)}$ have the same singular values.
}\label{K3pic}
\end{figure}

Figure~\ref{K3pic} shows the isochrone lines for the various
algorithms applied to the operator $K^{(3)}$. To make
comparison with figure~\ref{isochronefig3} straightforward, the
total time span is again $1$ minute subdivided in $6$s
intervals. Convergence first progresses faster than in
figure~\ref{isochronefig3}: the isochrone corresponding to
$t=6$s lies lower than in figure~\ref{isochronefig3}. For
$\lambda<\lambda_\mathrm{max}/2^6$ it is clear that the various
algorithms perform worse for the operator $K^{(3)}$ than for
the matrix $K^{(2)}$. As these two operators have identical
spectra, this implies that a well-conditioned spectrum alone is
no guarantee for good performance of the algorithms.

Figure~\ref{K3pic} shows the isochrone plots for the operator
$K^{(4)}$ which has identical spectrum as $K^{(1)}$, implying
that it is very ill-conditioned. The singular vectors of
$K^{(4)}$ are the same as for the random Gaussian matrix
$K^{(2)}$. We see that the ill-conditioning of the spectrum
influences the performance of the algorithms is a negative way,
as compared to the Gaussian operator $K^{(2)}$.

\begin{figure}
\centering \resizebox{\textwidth}{!}{\includegraphics{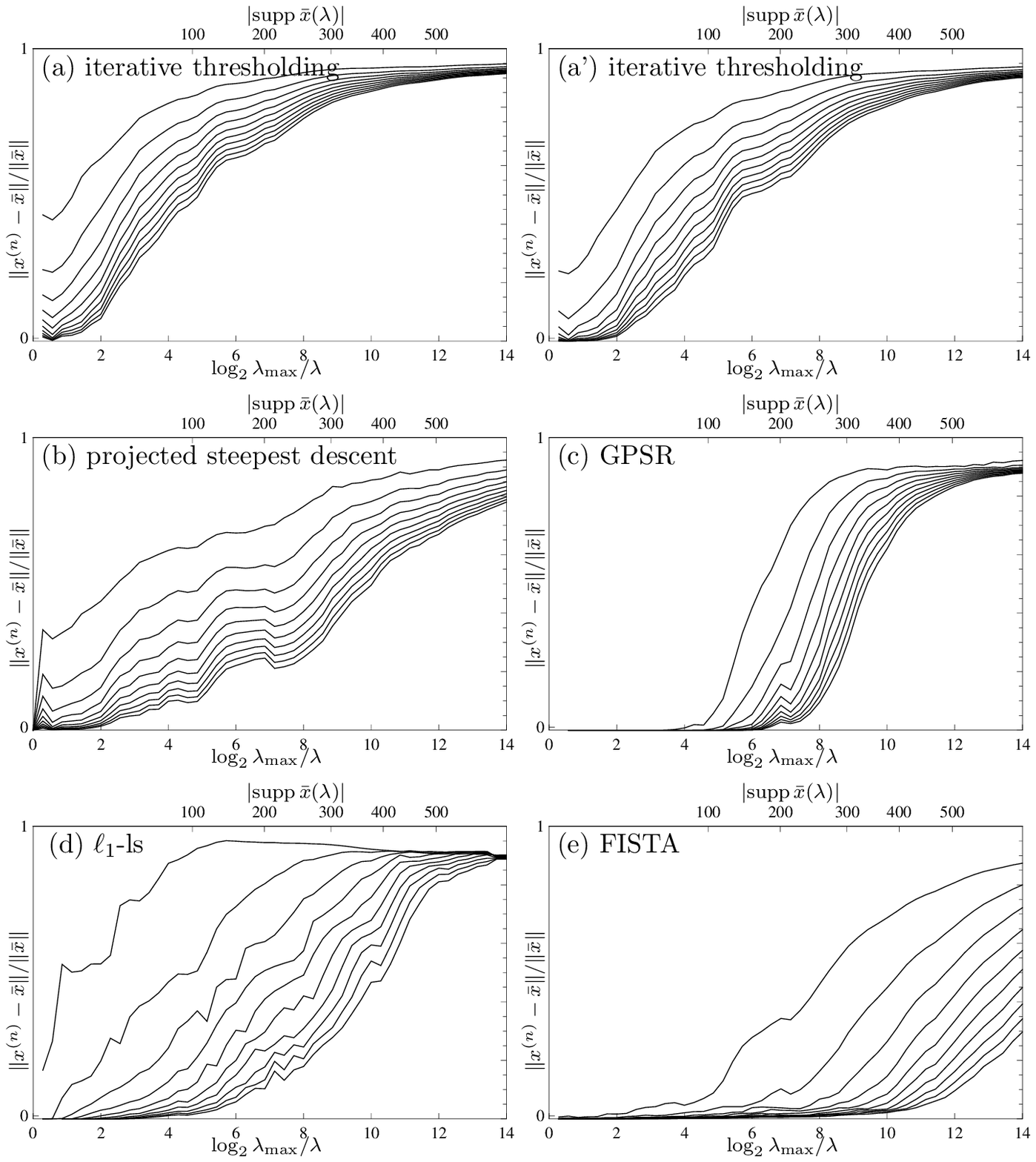}}
\caption{The six algorithms (a)--(e) are compared in case of the operator $K^{(4)}$.
The time intervals are $1,2,\ldots,10$m as opposed to $6,\ldots,60$s in figure~\ref{isochronefig3}. Clearly, replacing the spectrum of the random
Gaussian matrix $K^{(2)}$ by the spectrum of $K^{(1)}$, and leaving the singular vectors unaltered,
changes the behavior of the minimization algorithms in a negative way.
Convergence is much slower than in figure~\ref{isochronefig3} for all algorithms.}\label{K4pic}
\end{figure}

\subsection{Warm-start strategies}

For the warm-start algorithms (A)--(B), it is not possible to
draw isochrones because these methods are not purely iterative.
They depend on a preset maximum number of iterates and a preset
end value for the penalty parameter $\lambda$ or the $\ell_1$
norm $\rho$ of $x$. It is possible, for a fixed total
computation time and a fixed value of $\lambda_\mathrm{stop}$
or $\rho_\mathrm{stop}$, to plot $\|\bar
x(\lambda_n)-x^{(n)}\|$ vs. $\lambda_n$ or $\|\tilde
x(\rho_n)-x^{(n)}\|$ vs. $\rho_n$. This gives a condensed
picture of the performance of such an algorithm, as it include
information on the remaining error for various values of
$\lambda$ or $\rho$.

\begin{figure}
\centering\resizebox{\textwidth}{!}{\includegraphics{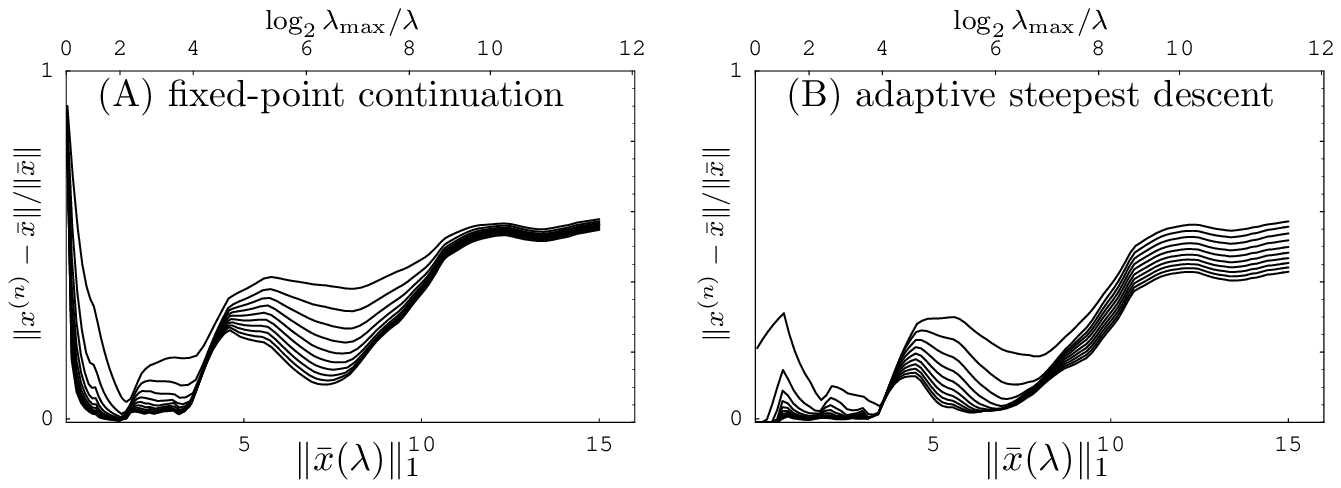}}
\caption{The behavior of two warm-start strategies for the operator $K^{(1)}$. The labels (A)--(B) refer to the list
in section \ref{problemsection}. In each picture the top line
corresponds to a total calculation time of $1$ minute. The
bottom line corresponds to $10$ minutes total calculation time
(this is about $4800$ iterative soft-thresholding steps or
$3000$ projected steepest descent steps). The lines in these
plots are not isochrones, as these algorithms are not purely
iterative, but depend on a preset stopping radius
$\rho_\mathrm{stop}$ or stopping penalty $\lambda_\mathrm{stop}$.
$\rho_\mathrm{stop}$ equals $15$ in this case.}\label{fpcfig}
\end{figure}

In figure \ref{fpcfig}, the warm-start methods (A) and (B) are
compared in the range $0\leq \|\bar x\|_1\leq 15$ for the
matrix $K^{(1)}$. For each experiment, ten runs were performed,
corresponding to total computation times of $1, 2, \ldots, 10$
minutes. For each run, the parameter $\rho_\mathrm{max}$ was
chosen to be 15 in algorithm (B). For algorithm (A),
$\lambda_\mathrm{stop}$ was chosen such that $\|\bar
x(\lambda_\mathrm{stop})\|_1=15$ also. From the pictures we see
that the algorithms (A)--(B) do not do very well for small
values of the penalty parameter $\lambda$. Their big advantage
is clear: For the price of a single run, an acceptable
approximation of \emph{all} minimizers $\bar x(\lambda)$ with
$\lambda_\mathrm{max}\geq\lambda\geq\lambda_\mathrm{stop}$ is
obtained. The algorithm (B) does somewhat better than algorithm
(A).

\begin{figure}
\centering\resizebox{\textwidth}{!}{\includegraphics{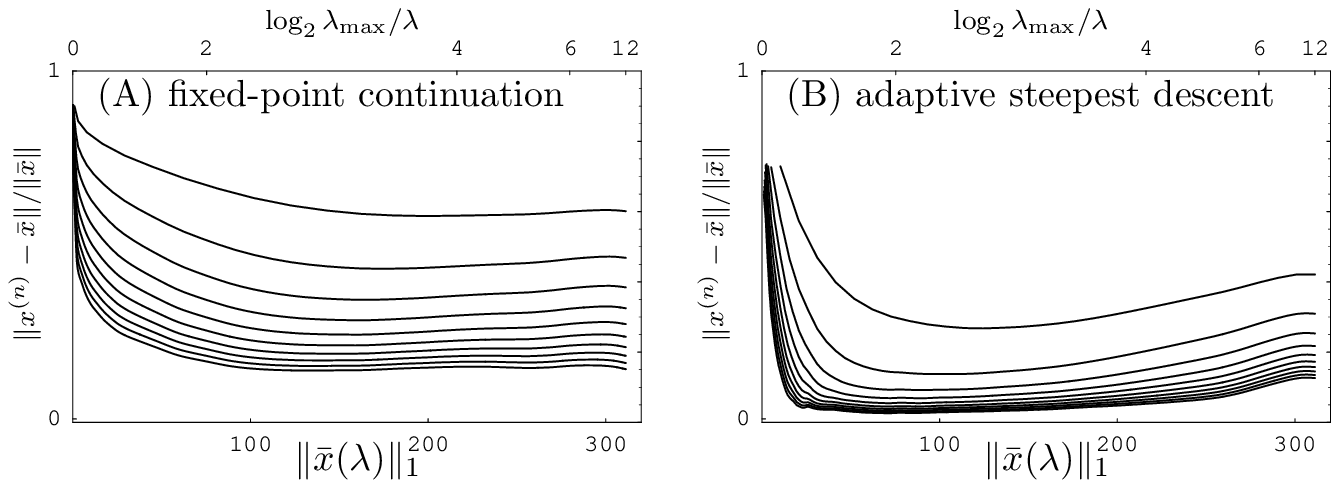}}
\caption{The behavior of two warm-start strategies for the
Gaussian random matrix $K^{(2)}$. In each picture the
top line corresponds to a total calculation time of $6$
seconds. The bottom line to $1$ minute. In this case we chose
to stop at $\lambda_\mathrm{stop}$ equals
$\lambda_\mathrm{max}/2^{12}$.}\label{gaussfpcfig}
\end{figure}

In figure \ref{gaussfpcfig}, the same type of comparison is
made for the matrix $K^{(2)}$. In this case, ten runs are
performed with a total computation time per run equal to
$6,12,\ldots,60$ seconds ($60$s corresponds to about $460$
iterative soft-thresholding steps or $300$ projected steepest
descent steps). This is ten times less than in case of the
matrix $K^{(1)}$. Both the `fixed point continuation method'
(A) and the `projected steepest descent' (B) do acceptably
well, even up to very small values for the penalty parameter
$\lambda$ (large value of $\rho$).

All the calculations in this note were done on a 2GHz processor and
2Gb ram memory.

\section{Conclusions}

The problem of assessing the convergence properties for $\ell_1$
penalized least squares functionals was discussed.

We start from the rather obvious observation that convergence
speed can only refer to the behavior of $e=\|x^{(n)}-\bar
x\|/|\bar x\|$ as a function of time. Luckily, in the case of
functional (\ref{l1functional}), the exact minimizer $\bar x$
(up to computer round-off) can be found in a finite number of
steps: even though the variational equations (\ref{kkt}) are
nonlinear, they can still be solved exactly using the
LARS/homotopy method. A direct calculation of the minimizers
$\bar x(\lambda)$ is thus possible, at least when the number of
nonzero coefficients in $\bar x$ is not too large (e.g. in the
numerical experiments in \cite{Figueiredo.Nowak.ea2008}, one
has $4096$ degrees of freedom and only about 160 nonzeros.).

We provided a graph that indicates the time complexity of the
exact algorithm as a function of $s$, the number of nonzero
components in $\bar x(\lambda)$. It showed us that computing
time rises approximately cubically as a function of $s$ (it is
linear at first). Also we gave an example where the size of
support of $\bar x(\lambda)$ does not decrease monotonically as
$\lambda$ decreases. The direct method is certainly practical
for $|\,\mathrm{supp}\, \bar x(\lambda)|\leq 10^3$ or so.

It is impossible to completely characterize the performance of
an iterative algorithm in just a single picture. A good
qualitative appreciation, however, can be gained from looking
at the approximation isochrones introduced in this note. These
lines in the $\lambda-e$-plane tell us for which value of the
penalty parameter $\lambda$ convergence is adequately fast, and
for which values it is inacceptably slow. We look at the region
$e\in [0, 1]$, because this is probably most interesting for
doing real inversions. One could also use a logarithmic scale
for $e$, and look at very small values of $e$ approaching
computer epsilon. But that is probably not of principle
interest to people doing real inverse problems. The main
content is thus in the concept of approximation isochrone and
in figures \ref{singleisochronefig}, \ref{isochronefig},
\ref{isochronefig3}, \ref{K3pic} and \ref{K4pic} comparing six
different algorithms for four operators.

For large penalty parameters, all algorithms mentioned in this
note do well for our particular example operators with a small
preference for the GPSR method. The biggest difference can be
found for small penalty parameters. Algorithms (a), (a'), (b)
and (c) risk to be useless in this case. The $\ell_1$-ls (d)
algorithm seems to be more robust, but it loses out for large
penalties. The FISTA method (e) appears to work best for small
values of the penalty parameter.

We uncovered two aspects that may influence the convergence
speed of the iterative algorithms. Firstly we saw that
convergence is slower for ill-conditioned operators than for
well-conditioned operators (comparison of $K^{(2)}$ and
$K^{(4)}$). Furthermore, the number of nonzero coefficients
that a recoverable minimizer has is smaller for an
ill-conditioned operator. In particular, the operator $K^{(1)}$
that comes out of a real physical problem, presents a challenge
for all algorithms for small values of the penalty parameter.
Secondly, the orientation of the null-space with respect to the
edges of the $\ell_1$ ball also influences the speed of
convergence of the iterative algorithms, even if the operator
is well-conditioned (comparison of $K^{(2)}$ and $K^{(3)}$).

We also compared two warm-start strategies and showed that
their main advantage is to yield a whole set of minimizers (for
different penalty parameters) in a single run. We found that
the adaptive-steepest descent method (B), proposed in
\cite{DaFoL2008} but tested here for the first time, does
better than the fixed point continuation method (A).

\section{Acknowledgments}

Part of this work was done as a post-doctoral research fellow
for the F.W.O-Vlaanderen (Belgium) at the Vrije Universiteit
Brussel and part was done as `Francqui Foundation
intercommunity post-doctoral researcher' at D\'epartement de
Math\'ematique, Universit\'e Libre de Bruxelles. Discussions
with Ingrid Daubechies and Christine De Mol are gratefully
acknowledged. The author acknowledges the financial support of
the VUB through the GOA-62 grant, of the FWO-Vlaanderen through
grant G.0564.09N.


\bibliographystyle{unsrt}
\bibliography{performance}{}

\begin{thebibliography}{10}

\bibitem{Donoho2006}
D.L. Donoho.
\newblock Compressed sensing.
\newblock {\em Information Theory, IEEE Transactions on}, 52(4):1289--1306,
  2006.

\bibitem{Daubechies.Defrise.ea2004}
I.~Daubechies, M.~Defrise, and C.~De~Mol.
\newblock An iterative thresholding algorithm for linear inverse problems with
  a sparsity constraint.
\newblock {\em Communications On Pure And Applied Mathematics},
  57(11):1413--1457, November 2004.

\bibitem{Santosa.Symes1986}
Fadil Santosa and William~W. Symes.
\newblock Linear inversion of band-limited reflection seismograms.
\newblock {\em SIAM J. Sci. Stat. Comput.}, 7:1307--1330, 1986.

\bibitem{Tibshirani1996}
R.~Tibshirani.
\newblock Regression shrinkage and selection via the lasso.
\newblock {\em J. Royal. Statist. Soc B.}, 58:267--288, 1996.

\bibitem{DaFoL2008}
I.~Daubechies, M.~Fornasier, and I.~Loris.
\newblock Accelerated projected gradient method for linear inverse problems
  with sparsity constraints.
\newblock {\em Journal of Fourier Analysis and Applications}, 2008.
\newblock Accepted.

\bibitem{Kim.Koh.ea2007}
S.-J. Kim, K.~Koh, M.~Lustig, S.~Boyd, and D.~Gorinevsky.
\newblock A method for large-scale {$\ell_1$}-regularized least squares
  problems with applications in signal processing and statistics.
\newblock {\em IEEE Journal on Selected Topics in Signal Processing}, 2007.
\newblock Accepted.

\bibitem{Figueiredo.Nowak.ea2008}
Mario A.~T. Figueiredo, Robert~D. Nowak, and Stephen~J. Wright.
\newblock Gradient projection for sparse reconstruction: Application to
  compressed sensing and other inverse problems.
\newblock {\em To appear in the IEEE Journal of Selected Topics in Signal
  Processing: Special Issue on Convex Optimization Methods for Signal
  Processing}, 2008.

\bibitem{Hale.Yin.ea2007}
Elaine~T. Hale, Wotao Yin, and Yin Zhang.
\newblock A fixed-point continuation method for {$\ell_1$}-regularized
  minimization with applications to compressed sensing.
\newblock Technical report, Rice University, 2007.

\bibitem{Beck.Teboulle2008}
Amir Beck and Marc Teboulle.
\newblock A fast iterative shrinkage-thresholding algorithm for linear inverse
  problems.
\newblock {\em SIAM Journal on Imaging Sciences}, 2008.
\newblock To appear.

\bibitem{Loris.Nolet.ea2007}
Ignace Loris, Guust Nolet, Ingrid Daubechies, and F.~A. Dahlen.
\newblock Tomographic inversion using $\ell_1$-norm regularization of wavelet
  coefficients.
\newblock {\em Geophysical Journal International}, 170(1):359--370, 2007.

\bibitem{Osborne.Presnell.ea2000}
M.~R. Osborne, B.~Presnell, and B.~A. Turlach.
\newblock A new approach to variable selection in least squares problems.
\newblock {\em IMA J. Numer. Anal.}, 20(3):389--403, July 2000.

\bibitem{Efron.Hastie.ea2004}
Bradley Efron, Trevor Hastie, Iain Johnstone, and Robert Tibshirani.
\newblock Least angle regression.
\newblock {\em Ann. Statist.}, 32(2):407--499, 2004.

\bibitem{IMM2005-03897}
K.~Sj{\"{o}}strand.
\newblock Matlab implementation of {LASSO,} {LARS,} the elastic net and {SPCA},
  jun 2005.
\newblock Version 2.0.

\bibitem{Donoho.Stodden.ea2007}
David Donoho, Victoria Stodden, and Yaakov Tsaig.
\newblock About {S}parse{L}ab, March 2007.

\bibitem{Loris2007}
Ignace Loris.
\newblock {L}1{P}ackv2: A {M}athematica package for minimizing an
  $\ell_1$-penalized functional.
\newblock {\em Computer Physics Communications} 179: 895--902, 2008.

\bibitem{Combettes.Wajs2005}
Patrick~L. Combettes and Valerie~R. Wajs.
\newblock Signal recovery by proximal forward-backward splitting.
\newblock {\em Multiscale Model. Simul.}, 4(4):1168--1200, January 2005.

\bibitem{Combettes1997}
P.L. Combettes.
\newblock Convex set theoretic image recovery by extrapolated iterations of
  parallel subgradient projections.
\newblock {\em Image Processing, IEEE Transactions on}, 6(4):493--506, 1997.

\end{thebibliography}

\end{document}